\documentclass{article}

\usepackage{proof}
\usepackage{latexsym}
\usepackage{amssymb}

\usepackage[all]{xy}

\newcommand{\blem}{\begin{lemma}}
\newcommand{\elem}{\end{lemma}}
\newcommand{\bth}{\begin{theorem}}
\newcommand{\ethm}{\end{theorem}}
\newcommand{\benu}{\begin{enumerate}}
\newcommand{\eenu}{\end{enumerate}}
\newcommand{\bdes}{\begin{description}}
\newcommand{\edes}{\end{description}}
\newcommand{\bdf}{\begin{definition}}
\newcommand{\edf}{\end{definition}}
\newcommand{\bcor}{\begin{cor}}
\newcommand{\ecor}{\end{cor}}
\newcommand{\bprp}{\begin{proposition}}
\newcommand{\eprp}{\end{proposition}}
\newcommand{\bmlem}{\begin{mlemma}}
\newcommand{\emlem}{\end{mlemma}}
\newcommand{\bclm}{\begin{claim}}
\newcommand{\eclm}{\end{claim}}
\newcommand{\bprf}{{\bf Proof}.\hspace{2mm}}

\newcommand{\eprf}{\hspace*{\fill} $\Box$}

\newcommand{\beqn}{\begin{equation}}
\newcommand{\eeqn}{\end{equation}}
\newcommand{\beqnarr}{\begin{eqnarray}}
\newcommand{\eeqnarr}{\end{eqnarray}}
\newcommand{\beqnarrs}{\begin{eqnarray*}}
\newcommand{\eeqnarrs}{\end{eqnarray*}}

\newcommand{\spand}{\,\&\,}

\newcommand{\bfG}{\mbox{\boldmath$G$} }

\newtheorem{theorem}{Theorem}[section]
\newtheorem{definition}[theorem]{Definition}
\newtheorem{proposition}[theorem]{Proposition}
\newtheorem{lemma}[theorem]{Lemma}
\newtheorem{cor}[theorem]{Corollary}
\newtheorem{mlemma}[theorem]{Main Lemma}
\newtheorem{claim}[theorem]{Claim}

\newcommand{\alp}{\alpha}

\newcommand{\veps}{\varepsilon}
\newcommand{\del}{\delta}
\newcommand{\Del}{\Delta}
\newcommand{\ome}{\omega}

\newcommand{\bet}{\beta}
\newcommand{\gam}{\gamma}
\newcommand{\Gam}{\Gamma}

\newcommand{\sig}{\sigma}
\newcommand{\Sig}{\Sigma}

\newcommand{\lam}{\lambda}
\newcommand{\Lam}{\Lambda}
\newcommand{\vphi}{\varphi}

\newcommand{\fal}{\forall}
\newcommand{\exi}{\exists}

\newcommand{\Rarw }{\Rightarrow}

\newcommand{\lrarw}{\leftrightarrow}
\newcommand{\Lrarw}{\Leftrightarrow}

\newcommand{\calt}{{\cal T}}

\newcommand{\la}{\langle}
\newcommand{\ra}{\rangle}

\newcommand{\msfiv}{\mbox{\hspace{5mm}}}

\newcommand{\brem}{\begin{remark}}
\newcommand{\erem}{\end{remark}}
\newcommand{\dg}{\mbox{{\rm dg}}}
\newcommand{\Natural}{\mathbb{N}}

\begin{document}

\title{Proof-theoretic strengths of the well-ordering principles
}

\author{Toshiyasu Arai
\thanks{
I'd like to thank A. Freund for pointing out a flaw in \cite{Arai17}
}
\\
Graduate School of Mathematical Sciences,
The University of Tokyo,
\\
3-8-1 Komaba, Meguro-ku,
Tokyo 153-8914, JAPAN
\\
tosarai@ms.u-tokyo.ac.jp}

\maketitle

\begin{abstract}
In this note the proof-theoretic ordinal of
the well-ordering principle for the normal functions ${\sf g}$ on ordinals
is shown to be equal to the least fixed point of ${\sf g}$.
Moreover corrections to the previous paper \cite{Arai17} are made.

\end{abstract}

\section{Introduction}

In this note
 we are concerned with a proof-theoretic strength of a $\Pi^{1}_{2}$-statement ${\rm WOP}({\sf g})$ 
saying that
`for any well-ordering $X$, ${\sf g}(X)$ is a well-ordering',
where ${\sf g}: \mathcal{P}(\Natural)\to\mathcal{P}(\Natural)$ is a computable functional on sets $X$ of natural numbers.
$\la n,m\ra$ denotes the elementary recursive pairing function 
$\la n,m\ra=\frac{(n+m)(n+m+1)}{2}+m$ on $\Natural$.

\bdf
{\rm
$X\subset\Natural$ defines a binary relation 
$<_{X}:=\{(n,m): \la n,m\ra\in X\}$.
\beqnarrs
{\rm Prg}[<_{X},Y] & :\Lrarw & \fal m\left(\fal n<_{X}m\,Y(n)\to Y(m)\right)
\\
{\rm TI}[<_{X},Y] & :\Lrarw & {\rm Prg}[<_{X},Y] \to\fal n\, Y(n)
\\
{\rm TI}[<_{X}] & :\Lrarw & \fal Y\,{\rm TI}[<_{X},Y]
\\
{\rm WO}(X) & :\Lrarw & {\rm LO}(X) \land {\rm TI}[<_{X}]
\eeqnarrs
where ${\rm LO}(X)$ denotes a $\Pi^{0}_{1}$-formula stating that $<_{X}$ is a linear ordering.

For a functional ${\sf g}: \mathcal{P}(\Natural)\to\mathcal{P}(\Natural)$,
\[
{\rm WOP}({\sf g}) :\Lrarw \fal X\left( {\rm WO}(X)\to {\rm WO}({\sf g}(X)) \right)
\]
}
\edf

The theorem due to J.-Y. Girard is a base for further results on the 
strengths of the well-ordering principles ${\rm WOP}({\sf g})$.
For second order arithmetics ${\rm RCA}_{0}$, ${\rm ACA}_{0}$, etc. see \cite{Simpson}.
For a set $X\subset\Natural$, $\ome^{X}$ denotes an ordering on $\Natural$ canonically defined
such that its order type is $\ome^{\alp}$ when $<_{X}$ is a well-ordering of type $\alp$.

\bth\label{th:Girardp.299}{\rm (Girard\cite{Girard})}
\\
Over ${\rm RCA}_{0}$, ${\rm ACA}_{0}$ is equivalent to ${\rm WOP}(\lam X.\ome^{X})$.
\end{theorem}

In the following theorem 
${\rm ACA}_{0}^{+}$ denotes an extension of ${\rm ACA}_{0}$ by the axiom of the existence of
the $\ome$-th jump of a given set.

\bth\label{th:MontalbanMarcone}{\rm (Marcone and Montalb\'an\cite{MontalbanMarcone})}
Over ${\rm RCA}_{0}$, ${\rm ACA}_{0}^{+}$ is equivalent to ${\rm WOP}(\lam X.\veps_{X})$.
\end{theorem}

Theorem \ref{th:MontalbanMarcone} is proved
in \cite{MontalbanMarcone} computability theoretically.
M. Rathjen noticed that the principle ${\rm WOP}({\sf g})$ is tied to the existence of
\textit{countable coded $\ome$-models}.

\bdf
{\rm
A \textit{countable coed $\ome$-model} of a second-order arithmetic $T$ is a set $Q\subset\Natural$ such that
$M(Q)\models T$, where
$M(Q)=\la\Natural,\{(Q)_{n}\}_{n\in\Natural}, +,\cdot,0,1,<\ra$
with $(Q)_{n}=\{m\in\Natural: \la n,m\ra\in Q\}$.

Let $X\in_{\ome}Y:\Lrarw(\exi n[X=(Y)_{n}])$ and
$X=_{\ome}Y:\Lrarw(\fal Z(Z\in_{\ome}X\lrarw Z\in_{\ome}Y))$.
}
\edf
It is not hard to see that over ${\rm ACA}_{0}$,
the existence of the $\ome$-th jump is equivalent to
the fact that there exists an arbitrarily large countable coded $\ome$-model of ${\rm ACA}_{0}$, cf.\,\cite{AfshariRathjen}.
The fact means that there is a countable coded $\ome$-model $Q$ of ${\rm ACA}_{0}$
containing a given set $X$, i.e., $X=(Q)_{0}$.
From this characterization, Afshari and Rathjen\cite{AfshariRathjen} gives 
a purely proof-theoretic proof of Theorem \ref{th:MontalbanMarcone}.
Their proof is based on Sch\"utte's method of complete proof search in $\ome$-logic, cf.\,\cite{Schutte}.

In \cite{MontalbanMarcone}, a further equivalence is established for the binary Veblen function.
In M. Rathjen, et.\,al.\cite{AfshariRathjen,RathjenWeiermannRM,RathjenHFriedman} and
\cite{Arai17} the well-ordering principles are investigated proof-theoretically.
Note that in Theorem \ref{th:Girardp.299} the proof-theoretic ordinal 
$|{\rm ACA}_{0}|=|{\rm WOP}(\lam X.\ome^{X})|=\veps_{0}$ is the least fixed point of the function 
$\lam x.\ome^{x}$.
Moreover the ordinal $|{\rm ACA}_{0}^{+}|=|{\rm WOP}(\lam X.\veps_{X})|$ in 
\cite{MontalbanMarcone,AfshariRathjen} is the least fixed point of the function $\lam x.\veps_{x}$, and
$|{\rm ATR}_{0}|=|{\rm WOP}(\lam X.\vphi X0)|=\Gam_{0}$ in \cite{RathjenWeiermannRM} 
one of $\lam x.\vphi_{x}(0)$.
These results suggest a general result that
the well-ordering principle for normal functions ${\sf g}$ on ordinals
is equal to the least fixed point of ${\sf g}$.

In this note we confirm this conjecture
under a mild condition on normal function ${\sf g}$, 
cf.\,Definition \ref{df:omegag} for the extendible term structures.

We assume that the strictly increasing function ${\sf g}$ enjoys the following conditions.
The computability of the functional ${\sf g}$ and the linearity of ${\sf g}(X)$ for linear orderings $X$ are assumed
to be provable elementarily,
and if $X$ is a well-ordering of type $\alp$, then ${\sf g}(X)$ is also a well-ordering of type ${\sf g}(\alp)$.
Moreover
${\sf g}(X)$ is assumed to be a \textit{term structure} over
constants ${\sf g}(c)\, (c\in X)$, function constants $+,\ome$, and possibly other function constants.

\bth\label{th:prfthordwop}
Let ${\sf g}(X)$ be an extendible term structure, and 
${\sf g}^{\prime}(X)$ an exponential term structure for which (\ref{df:omegag.4}) holds below.

Then the proof-theoretic ordinal of the second order arithmetic
${\rm WOP}({\sf g})$ over ${\rm ACA}_{0}$
is equal to the least fixed point ${\sf g}^{\prime}(0)$ of the ${\sf g}$-function,
$|{\rm ACA}_{0}+{\rm WOP}({\sf g})|=\min\{\alp:{\sf g}(\alp)=\alp\}=\min\{\alp>0:\fal\bet<\alp({\sf g}(\bet)<\alp)\}$.
\end{theorem}

On the other side the proof of the harder direction of Theorem 4 in \cite{Arai17} 
should be corrected as pointed out by A. Freund.
The theorem is stated as follows.

\bth\label{th:derivativemodel}
Let ${\sf g}(X)$ be an extendible term structure, and 
${\sf g}^{\prime}(X)$ an exponential term structure for which (\ref{df:omegag.4}) holds.

Then the following two are mutually equivalent over ${\rm ACA}_{0}$:
\benu
\item\label{th:derivativemodel.1}
${\rm WOP}({\sf g}^{\prime})$.
\item\label{th:derivativemodel.2}
$\left({\rm WOP}({\sf g})\right)^{+}
:\Lrarw \fal X\exi \mathcal{Q}[X\in_{\ome}\mathcal{Q} \land 
M(\mathcal{Q})\models{\rm ACA}_{0}+{\rm WOP}({\sf g})]$.
\eenu
\end{theorem}

Let us mention the contents of the paper.
In the next section \ref{sect:termstr}, ${\sf g}(X)$ is defined as a term structure.
Exponential term structures and extendible ones are defined.
The easy direction in Theorem \ref{th:prfthordwop} is shown.
In section \ref{sect:proofschema} we prove Theorems \ref{th:prfthordwop} and \ref{th:derivativemodel},
assuming an elimination theorem \ref{lem:WPbndtransf} of the well-ordering principle
in infinitary sequent calculi. 
In section \ref{sect:CE} we prove the elimination theorem \ref{lem:WPbndtransf}.

\section{Term structures}\label{sect:termstr}
Let us reproduce definitions on term structure from \cite{Arai17}.

The fact that ${\sf g}$ sends
linear orderings $X$ to linear orderings ${\sf g}(X)$ should be provable in an elementary way.
${\sf g}$ sends a binary relation $<_{X}$ on a set $X$ to a binary relation
$<_{{\sf g}(X)}={\sf g}(<_{X})$ on a set ${\sf g}(X)$.
We further assume that ${\sf g}(X)$ is a Skolem hull, i.e., a term structure over
constants $0$ and ${\sf g}(c)\, (c\in \{0\}\cup X)$ with the least element $0$ in the order $<_{X}$,
the addition $+$, the exponentiation $\ome^{x}$, and possibly other function constants 
in a list $\mathcal{F}$.
When $\mathcal{F}=\emptyset$, let $\ome^{\alp}:={\sf g}(\alp)$.
Otherwise we assume that $\lam \xi.\,\ome^{\xi}$ is in the list $\mathcal{F}$.

\bdf\label{df:gtmstr}
{\rm
\benu
\item
${\sf g}(X)$ is said to be a \textit{computably linear} term structure if
there are three $\Sig^{0}_{1}(X)$-formulas ${\sf g}(X), <_{{\sf g}(X)},=$ 
for which all of the following facts are provable in ${\rm RCA}_{0}$:
let $\alp,\bet,\gam,\ldots$ range over terms.
\benu
\item(Computability)
Each of ${\sf g}(X)$, $<_{{\sf g}(X)}$ and $=$ is $\Del^{0}_{1}(X)$-definable.
${\sf g}(X)$ is a computable set, and $<_{{\sf g}(X)}$ and $=$ are computable binary relations.

\item(Congruence)
\\
$=$ is a congruence relation on the structure $\la{\sf g}(X);<_{{\sf g}(X)},f,\ldots\ra$.

Let us denote ${\sf g}(X)/=$ the quotient set.

In what follows assume that $<_{X}$ is a linear ordering on $X$.

 \item(Linearity)
 $<_{{\sf g}(X)}$ is a linear ordering on ${\sf g}(X)/=$ with the least element $0$.

\item(Increasing)
${\sf g}$ is strictly increasing:
$c<_{X}d \Rarw {\sf g}(c)<_{{\sf g}(X)}{\sf g}(d)$.

\item(Continuity)
${\sf g}$ is continuous:
Let $\alp<_{{\sf g}(X)}{\sf g}(c)$ for a limit $c\in X$ and $\alp\in{\sf g}(X)$.
Then there exists a $d<_{X}c$ such that $\alp<_{{\sf g}(X)}{\sf g}(d)$.
\eenu

\item
A computably linear term structure ${\sf g}(X)$ is said to be \textit{extendible} if it enjoys
the following two conditions.
\benu
\item(Suborder)
If $\la X,<_{X}\ra$ is a substructure of $\la Y,<_{Y}\ra$, then 
$\la{\sf g}(X);=, <_{{\sf g}(X)}, f,\ldots\ra$ is a substructure of $\la{\sf g}(Y); =, <_{{\sf g}(Y)},f,\ldots\ra$.

\item(Indiscernible)
\\
$\la {\sf g}(c): c\in\{0\}\cup X\ra$ is an indiscernible sequence for 
linear orderings $\la{\sf g}(X),<_{{\sf g}(X)}\ra$:
Let $\alp[0,{\sf g}(c_{1}),\ldots,{\sf g}(c_{n})], \bet[0,{\sf g}(c_{1}),\ldots,{\sf g}(c_{n})]\in{\sf g}(X)$ be
terms such that constants occurring in them are among the list $0,{\sf g}(c_{1}),\ldots,{\sf g}(c_{n})$.
Then for any increasing sequences $c_{1}<_{X}\ldots<_{X}c_{n}$ and
$d_{1}<_{X}\ldots<_{X}d_{n}$, the following holds.
\beqnarr
&&
\alp[0,{\sf g}(c_{1}),\ldots,{\sf g}(c_{n})]<_{{\sf g}(X)}\bet[0,{\sf g}(c_{1}),\ldots,{\sf g}(c_{n})]
\label{eq:indiscerng}
\\
& \Lrarw &
\alp[0,{\sf g}(d_{1}),\ldots,{\sf g}(d_{n})]<_{{\sf g}(X)}\bet[0,{\sf g}(d_{1}),\ldots,{\sf g}(d_{n})]
\nonumber
\eeqnarr

\eenu
\eenu
}
\edf

\bprp\label{prp:dilate}
{\rm Suppose ${\sf g}(X)$ is an extendible term structure.
Then the following is provable in ${\rm RCA}_{0}$:}
Let both $X$ and $Y$ be linear orderings.

Let $f:\{0\}\cup X\to \{0\}\cup Y$ be an order preserving map, $n<_{X}m\Rarw f(n)<_{Y}f(m)\, (n,m\in\{0\}\cup X)$.
Then there is an order preserving map $F:{\sf g}(X)\to{\sf g}(Y)$, $n<_{{\sf g}(X)}m\Rarw F(n)<_{{\sf g}(Y)}F(m)$,
which extends $f$ in the sense that $F({\sf g}(n))={\sf g}(f(n))$.
\eprp
\bprf
This is seen from the indiscernibility (\ref{eq:indiscerng}), cf.\,\cite{Arai17}.
\eprf

\bdf\label{df:omegag}
{\rm Suppose that function symbols $+, \lam \xi.\,\ome^{\xi}$ are in the list $\mathcal{F}$ of function symbols for
a computably linear term structure ${\sf g}(X)$. Let $1:=\ome^{0}$, and $2:=1+1$, etc.

${\sf g}(X)$ is said to be an \textit{exponential} term structure 
(with respect to function symbols $+,\lam \xi.\,\ome^{\xi}$) if all of the followings are provable in ${\rm RCA}_{0}$.
\benu
 \item\label{df:omegag.1}
 $0$ is the least element in $<_{{\sf g}(X)}$, and $\alp+1$ is the successor of $\alp$.
 \item\label{df:omegag.2}
 $+$ and $\lam \xi.\,\ome^{\xi}$ enjoy the following familiar conditions.
 \benu
 \item
  $\alp<_{{\sf g}(X)}\bet\to \ome^{\alp}+\ome^{\bet}=\ome^{\bet}$.
  \item
 $\gam+\lam=\sup\{\gam+\bet:\bet<\lam\}$ when $\lam$ is a limit number, i.e.,
 $\lam\neq 0$ and $\fal \bet<_{{\sf g}(X)}\lam(\bet+1<_{{\sf g}(X)}\lam)$.
 \item
  $\bet_{1}<_{{\sf g}(X)}\bet_{2} \to \alp+\bet_{1}<_{{\sf g}(X)}\alp+\bet_{2}$, and
  $\alp_{1}<_{{\sf g}(X)}\alp_{2} \to \alp_{1}+\bet\leq_{{\sf g}(X)}\alp_{2}+\bet$.
  \item
  $(\alp+\bet)+\gam=\alp+(\bet+\gam)$.
  \item
  $\alp<_{{\sf g}(X)}\bet\to\exi\gam\leq_{{\sf g}(X)}\bet(\alp+\gam=\bet)$.
  \item
  Let $\alp_{n}\leq_{{\sf g}(X)}\cdots\leq_{{\sf g}(X)}\alp_{0}$ and 
  $\bet_{m}\leq_{{\sf g}(X)}\cdots\leq_{{\sf g}(X)}\bet_{0}$.
  Then $\ome^{\alp_{0}}+\cdots+\ome^{\alp_{n}}<_{{\sf g}(X)}\ome^{\bet_{0}}+\cdots+\ome^{\bet_{m}}$ iff
  either $n<m$ and $\fal i\leq n(\alp_{i}=\bet_{i})$, or 
  $\exi j\leq\min\{n,m\}[\alp_{j}<_{{\sf g}(X)}\bet_{j}\land \fal i<j(\alp_{i}=\bet_{i})]$.
 \eenu

\item\label{df:omegag.3}
Each $f(\bet_{1},\ldots,\bet_{n})\in{\sf g}(X)\, (+\neq f\in\mathcal{F})$ as well as 
${\sf g}(c)\,(c\in\{0\}\cup X)$ 
is closed under $+$.
In other words the terms $f(\bet_{1},\ldots,\bet_{n})$ and ${\sf g}(c)$ denote additively closed ordinals
 (additive principal numbers)
when $<_{{\sf g}(X)}$ is a well-ordering.

\eenu
}
\edf

In what follows we assume that ${\sf g}(X)$ is an extendible term structure,
and ${\sf g}^{\prime}(X)$ is an exponential term structure.
Constants in the term structure ${\sf g}^{\prime}(X)$ are $0$ and
${\sf g}^{\prime}(c)$ for $c\in\{0\}\cup X$,
and function symbols in $\mathcal{F}\cup\{0,+\}\cup\{{\sf g}\}$ with a unary function symbol ${\sf g}$.
We are assuming that a function constant $\lam \xi.\,\ome^{\xi}$ is in the list $\mathcal{F}\cup\{{\sf g}\}$.
Furthermore assume that ${\rm RCA}_{0}$ proves that
\beqnarr
\bet_{1},\ldots,\bet_{n}<_{{\sf g}^{\prime}(X)}{\sf g}^{\prime}(c) & \to &
 f(\bet_{1},\ldots,\bet_{n})<_{{\sf g}^{\prime}(X)}{\sf g}^{\prime}(c) \,(f\in\mathcal{F}\cup\{+,{\sf g}\})
 \nonumber
 \\
\ome^{{\sf g}^{\prime}(\bet)} & = & {\sf g}({\sf g}^{\prime}(\bet))={\sf g}^{\prime}(\bet)
\nonumber
\\
{\sf g}^{\prime}(0) & = & \sup_{n}{\sf g}^{n}(0)
\label{df:omegag.4}
\\
{\sf g}^{\prime}(c+1) & = & \sup_{n}{\sf g}^{n}({\sf g}^{\prime}(c)+1)\, (c\in\{0\}\cup X)
\nonumber
\eeqnarr
where ${\sf g}^{n}$ denotes the $n$-th iterate of the function ${\sf g}$,
and we are assuming in the last  that
the successor element $c+1$ of $c$ in $X$ exists.
The last two in (\ref{df:omegag.4}) hold for normal functions ${\sf g}$ when ${\sf g}(0)>0$.

Note that ${\sf g}^{\prime}(c)$ is an epsilon number when $<_{{\sf g}^{\prime}(X)}$
is a well-ordering since the exponential function is in $\mathcal{F}\cup\{{\sf g}\}$.
\\

We show the easy direction in Theorem \ref{th:prfthordwop}.
Let $<$ be an order of type ${\sf g}^{\prime}(0)$, which is defined from a family of structures
${\sf g}(X_{n})$ where the order types of $X_{n}$ is $\gam_{n}+1$ defined as follows.
A series of ordinals $\{\gam_{n}\}_{n}<{\sf g}^{\prime}(0)$ is defined recursively by $\gam_{0}=0$ and 
$\gam_{n+1}={\sf g}(\gam_{n})$.
Then ${\rm WOP}({\sf g})$ yields inductively ${\rm TI}[<_{\gam_{n}}]$ for initial segments of type $\gam_{n}$.
Hence $|{\rm WOP}({\sf g})|\geq{\sf g}^{\prime}(0):=\min\{\alp>0:\fal\bet<\alp({\sf g}(\bet)<\alp)\}$.

\section{Proof schema}\label{sect:proofschema}
In this section we give a proof schema of Theorems \ref{th:prfthordwop} and \ref{th:derivativemodel},
each of these is based on an elimination theorem \ref{lem:WPbndtransf} of the well-ordering principle
in infinitary sequent calculi. 


Formulas in our infinitary sequent calculi are generated from
literals $\top$(truth), $\bot:\equiv\bar{\top}$(absurdity),
$P(n),\bar{P}(n),E_{i}(n),\bar{E}_{i}(n), X_{i}(n),\bar{X}_{i}(n)\,(i,n\in\Natural)$
by applyig
infinitary disjunction $\bigvee_{n\in\Natural}A_{n}$, 
infinitary conjunction $\bigwedge_{n\in\Natural}A_{n}$
and second-order quantifications $\exi X,\fal X$.
Binary disjunctions $A_{0}\lor A_{1}$ are understood to be
$\bigvee_{n}B_{n}$ with $B_{0}\equiv A_{0}$ and $B_{1+n}\equiv A_{1}$,
and similarly for binary conjunctions.
A formula is said to be a well-formed formula, wff in short if
there is no free occurrence of `bound variables' $X_{i},\bar{X}_{i}$ in it.
 The negation $A$ of a wff $A$ is defined recursively by the de Morgan's law and the elimination
 of double negations.
Each wff is assumed to be a translation $A^{\infty}$ of
a formula $A$ without free first-order variables in the language of second-order arithmetic.
The translation is defined recursively as follows.
For an arithmetic literal $L$,
 $L^{\infty}\equiv \top$ if $L$ is true in the standard model $\Natural$,
 $L^{\infty}\equiv\bot$ otherwise.
 For a closed terms $t$ and $R\in\{P,\bar{P},E_{i},\bar{E}_{i},X_{i},\bar{X}_{i}: i\in\Natural\}$,
 $R(t)^{\infty}\equiv R(n)$ with the value $n$ of the closed term $t$ in $\Natural$.
 $(A_{0}\lor A_{1})^{\infty}\equiv(A_{0}^{\infty}\lor A_{1}^{\infty})$, and similarly for conjunctions.
 $(\exi x\,A(x))^{\infty}\equiv\bigvee_{n}A(\bar{n})^{\infty}$ for the $n$-th numeral $\bar{n}$.
 $(\fal x\,A(x))^{\infty}$ is defined to be an infinitary conjunction similarly.
 $(\exi X\,A(X))^{\infty}\equiv(\exi X\, A(X)^{\infty})$, and similarly for 
 the second-order universal quantifiers.
A formula is said to be a \textit{first-order} if no second-order quantifier occurs in it,
while it is \textit{arithmetical} if it is the translation of 
a formula in the language of the first-order arithmetic.
i.e., neither the predicate constant $P$ nor second-order variable occurs in it.

Each first-order formula $A$ defines a binary relation $n<_{A}m:\Lrarw A(\la n,m\ra)$.
The principle is formulated in the inference rule $(WP)$ together with a rule
for the progressiveness ${\rm Prg}[<_{A},E_{A}]$ of $E_{A}$ with respect to $<_{A}$:
\[
\infer[(WP)]{\Gam}
{
\{
\Gam, E_{A}(n)
: n\in\Natural
\}
&
\lnot{\rm TI}[<_{{\sf g}_{A}}],\Gam
}
\]
where
$E_{A}$ is a variable proper to the relation $<_{A}$, and does not occur in $\Gam$.
$n<_{{\sf g}_{A}}m:\Lrarw{\sf g}(A)(\la n,m\ra)$.

Our proof proceeds as follows.
Given cut-free derivations of $\Gam, E_{A}(n)$ without the rule $(WP)$,
suppose that we can obtain an embedding $f$ from the relation $<_{A}$ to an ordinal $\alp$
such that $n<_{A}m\Rarw f(n)<f(m)<\alp$.
Then the embedding $f$ can be extended to an embedding $F$ from the relation
$<_{{\sf g}_{A}}$ to an ordinal ${\sf g}(\alp)$ by Proposition \ref{prp:dilate}.
The embedding $F$ yields the transfinite induction ${\rm TI}[<_{{\sf g}_{A}}]$ 
for the relation $<_{{\sf g}_{A}}$.
Eliminating the false formula $\lnot{\rm TI}[<_{{\sf g}_{A}}]$,
we obtain $\Gam$.

However in order to extract such an embedding $f$ from
derivations,
we have to fix a meaning of the relation $<_{A}$.
In other words, we need to interpret the predicate constant $P$ and free-variables $E_{i}$
occurring in the formula $A$ so that these denote sets of natural numbers.
This motivates Definition \ref{df:GQ} below.


\bdf\label{df:GQ}
{\rm
Let $\mathcal{E}\subset\Natural$ be a family of sets
$\mathcal{E}_{i}=\{n\in\Natural: \la i,n\ra\in\mathcal{E}\}$.
Each variable $E_{i}$ is understood to denote 
the set $\mathcal{E}_{i}$.
Let
\[
{\rm Diag}(\mathcal{E}_{i})=
\{E_{i}(n): n\in\mathcal{E}_{i}\}\cup
\{\bar{E}_{i}(n): n\not\in\mathcal{E}_{i}\}.
\]

The predicate $P$
denotes a set $\mathcal{P}\subset\Natural$.
\[
{\rm Diag}(\mathcal{P})=\{P(n): n\in\mathcal{P}\}\cup\{\bar{P}(n) : n\not\in\mathcal{P}\}.
\]

${\rm Diag}(\mathcal{P},\mathcal{E})={\rm Diag}(\mathcal{P})\cup\bigcup_{i\in\Natural}{\rm Diag}(\mathcal{E}_{i})$ is identified with the countable coded $\ome$-model 
$\la\Natural;\mathcal{P},\mathcal{E}_{i}\ra_{i\in\Natural}$, and
${\rm Diag}(\mathcal{P},\mathcal{E})\models A:\Lrarw\la\Natural;\mathcal{P},\mathcal{E}\ra_{i\in\Natural}\models A$ for 
first-order formulas $A$.
For $\Sig^{1}_{1}$-formulas $\exi X\, F(X)$ with first-order matrices $F$,
define
${\rm Diag}(\mathcal{P},\mathcal{E})\models\exi X\, F(X)$ iff
there exists a first-order formula $A(x)$ in the language of arithmetic
such that
${\rm Diag}(\mathcal{P},\mathcal{E})\models F(A)$,
where
$F(A)$ denotes the result of replacing literals $X(n)$ [$\bar{X}(n)$] in $F(X)$ by
$A(n)^{\infty}$ [by $\lnot A(n)^{\infty}$], resp.

For a finite set $\Gam$ of first-order formulas, 
$Var(\Gam)$ denotes the set of second-order variables $E_{i}$
occurring in $\Gam$.
For a family $\mathcal{E}^{X}$ of finite sets $\mathcal{E}^{X}_{i}$, let
\[
\Del(\mathcal{E}^{X};\Gam):= \{E_{i}(n): n\in\mathcal{E}^{X}_{i}, E_{i}\in Var(\Gam),i\in\Natural\}
\]
}
\edf

\bdf\label{df:calculi}
{\rm
Let $\mathcal{P}\subset\Natural$ be a set of natural numbers, and
 $\mathcal{E}$ a family of sets $\mathcal{E}_{i}\subset\Natural$.
We define two cut-free infinitary one-sided sequent calculi
${\rm Diag}(\mathcal{P})+(prg)^{\infty}+(WP)$, and
${\rm Diag}(\mathcal{P},\emptyset)+(prg)^{\emptyset}$ as follows.

Let 
$\mathcal{E}^{X}$ be a family of \textit{finite} sets 
$\mathcal{E}^{X}_{i}\subset\Natural\,(i\in\Natural)$,
$\bet,\alp$ ordinals, and $\Gam$ a sequent, i.e., a finite set of formulas (in negation normal from).
We define a derivability relation
$\vdash^{\bet}_{\alp}\Gam$ in the calculus
${\rm Diag}(\mathcal{P})+(prg)^{\infty}+(WP)$,
and one
$\mathcal{E}^{X}\vdash^{\bet}\Gam$
in 
${\rm Diag}(\mathcal{P},\emptyset)+(prg)^{\emptyset}$
as follows,
where the depth of the derivation is bounded by $\bet$,
and the depth of the \textit{nested applications} of the inferences
 $(WP)$
is bounded by $\alp$ in the witnessed derivation in ${\rm Diag}(\mathcal{P})+(prg)^{\infty}+(WP)$.
\\
{\bf Axioms} or initial sequents:
\benu
\item
For $L\equiv\top$ and $L\in {\rm Diag}(\mathcal{P})$,
both $\vdash^{\bet}_{\alp}\Del,L$ 
and $\mathcal{E}^{X}\vdash^{\bet}\Del,L$ hold.

\item
$\vdash^{\bet}_{\alp}\Del,\bar{L},L$ for literals
$L\in\{E_{i}(n): i,n\in\Natural\}$.

\item
$\mathcal{E}^{X}\vdash^{\bet}\Del,\bar{L}$ for literals
$L\in\{E_{i}(n): i,n\in\Natural\}$.

\eenu
{\bf Inference rules}:  
The following inference rules $(\bigvee), (\bigwedge),(Rep),(\exi^{2}_{1st}),(\fal^{2})$
are shared by two calculi.
The left part $\mathcal{E}^{X}$ of $\vdash$ should be deleted for
the calculus ${\rm Diag}(\mathcal{P})+(prg)^{\infty}+(WP)$,
and the subscript $\alp$ is irrelevant to the calculus 
${\rm Diag}(\mathcal{P},\emptyset)+(prg)^{\infty}$ in the following.
Let $\gam<\bet$
\[
\infer[(\bigvee)]{\mathcal{E}^{X}\vdash^{\bet}_{\alp}\Gam,\bigvee_{n}A_{n}}
{\mathcal{E}^{X}\vdash^{\gam}_{\alp}\Gam,\bigvee_{n}A_{n},A_{i}}
\,\,
\infer[(\bigwedge)]{\mathcal{E}^{X}\vdash^{\bet}_{\alp}\Gam,\bigwedge_{n}A_{n}}
{
\{
\mathcal{E}^{X}\vdash^{\gam}_{\alp}\Gam,\bigwedge_{n}A_{n}, A_{i}
:i\in\Natural\}
}
\,\,
\infer[(Rep)]{\mathcal{E}^{X}\vdash^{\bet}_{\alp}\Gam}
{\mathcal{E}^{X}\vdash^{\gam}_{\alp}\Gam}
\]
\[
\infer[(\exi^{2}_{1st})]{\mathcal{E}^{X}\vdash^{\bet}_{\alp}\exi X F(X),\Gam}
{
\mathcal{E}^{X}\vdash^{\gam}_{\alp}F(A),\exi X F(X),\Gam}
\msfiv
\infer[(\fal^{2})]{\mathcal{E}^{X}\vdash^{\bet}_{\alp}\Gam,\fal X\, F(X)}
{
\mathcal{E}^{X}\vdash^{\gam}_{\alp}\Gam,\fal X\, F(X), F(E)
}
\]
where in $(\exi^{2}_{1st})$, $A(x)$ is a first-order formula, and
in $(\fal^{2})$, $E$ is an eigenvariable not occurring in $\Gam\cup\{\fal X F(X)\}$.

A first-order formula $A$ defines a binary relation $n<_{A}m:\Lrarw A(\la n,m\ra)$.
Let 
$n<_{{\sf g}_{A}}m:\Lrarw {\sf g}(A)(\la n,m\ra)$.
For each first-order formulas $A$, $\bet_{0}<\bet$ and $\alp_{0}<\alp$, we have
the following:
\[
\infer[(WP)]{\vdash^{\bet}_{\alp}\Gam}
{
\{
\vdash^{\bet_{0}}_{\alp_{0}}\Gam, E_{A}(n)
: n\in\Natural
\}
&
\vdash^{\bet_{0}}_{\alp_{0}}\lnot{\rm TI}[<_{{\sf g}_{A}}],\Gam
}
\]
where
the variable $E_{A}$ with the G\"odel number $i=\lceil A\rceil$ does not occur in $\Gam$.

\benu
\item
For each first-order formulas $A$, $\bet_{0}<\bet$, we have
the following.
Let $<_{A}^{*}$ denote the transitive closure of the relation $<_{A}$ for a first-order formula $A$.

\[
\infer[(prg)^{\infty}]{\vdash^{\bet}_{\alp}\Gam,E_{A}(m)}
{
\{
\vdash^{\bet_{0}}_{\alp}
\Gam,E_{A}(m),n\not<^{*}_{A}m,E_{A}(n)
:n\in\Natural
\}
}
\]
where $E_{A}\equiv E_{i}$ with the G\"odel number $i=\lceil A\rceil$ of the formula $A$,
and  $Var(A)\subset Var(\Gam)$.

\item
Let $A$ be a first-order formula with the G\"odel number $i=\lceil A\rceil$.
$n<^{*,\emptyset}_{A}m$ denotes the transitive closure of the relation 
$n<^{\emptyset}_{A}m:\Lrarw {\rm Diag}(\mathcal{P},\emptyset)\models A(\la n,m\ra)$.
If
\beqn\label{eq:normal}
m\in\mathcal{E}^{X}_{i}
\eeqn
then the inference $(prg)^{\emptyset}$ can be applied for $\bet_{0}<\bet$:
\[
\infer[(prg)^{\emptyset}]{((\mathcal{E}^{X}_{j})_{j\neq i},\mathcal{E}^{X}_{i})\vdash^{\bet}\Gam,E_{A}(m)}
{
\{((\mathcal{E}^{X}_{j})_{j\neq i},\mathcal{E}^{X}_{i}\cup\{n\})\vdash^{\bet_{0}}\Gam,E_{A}(m),E_{A}(n):n<^{*,\emptyset}_{A}m
\}
}
\]
where  
we assume that the variable $E_{A}\equiv E_{i}$ does not occur in $A$,
and $Var(A)\subset Var(\Gam)$.

\eenu
}
\edf
The rule $(prg)^{\infty}$ states the fact that the set
$E_{A}$ is progressive with respect to the relation 
$<^{*}_{A}$, i.e.,
${\rm Prg}[<^{*}_{A},E_{A}]$.
It is convenient for us
in proving Theorem \ref{th:TakeutimultiTI} in section \ref{sect:CE}
to have the weaker statement 
${\rm Prg}[<^{*}_{A},E_{A}]$ instead of the
stronger ${\rm Prg}[<_{A},E_{A}]$.
$(WP)$ together with $(prg)$
yields the well-ordering principle for ${\sf g}$.

\bdf\label{df:attach}
{\rm
Let $\pi$ be a derivation witnessing the fact
${\rm Diag}(\mathcal{P})+(prg)^{\infty}+(WP)\vdash^{\bet}_{\alp}\Gam_{0}$,
and $T(\pi)\subset{}^{<\ome}\Natural$ the underlying tree of $\pi$.
Let us assign recursively a family $\mathcal{E}^{X}(\sig)=\mathcal{E}^{X}(\sig;\pi)$
of \textit{finite} sets
$\mathcal{E}^{X}_{i}(\sig)$
to each node $\sig\in T(\pi)$ in a bottom-up way as follows.
In the definition,
$\mathcal{E}^{X}(\sig)\vdash\Gam$ designates that
$\mathcal{E}^{X}(\sig)$ is assigned to the node $\sig$ at which
the sequent $\Gam$ is placed.
In this case we write $\sig:\Gam$.

To the end-sequent, i.e., the bottom sequent $\emptyset:\Gam_{0}$,
assign the set $\mathcal{E}_{i}^{X}(\emptyset)=\{n: E_{i}(n)\in\Gam_{0}\}$.

Suppose that finite sets $(\mathcal{E}^{X}_{j}(\sig))_{j\neq i}$ are assigned to the lower sequent $\sig:\Gam$
of a rule $(WP)$ for the relation $<_{A}$ with $i=\lceil A\rceil$.
For the $n$-th left upper sequents $\sig*(n):\Gam,E_{A}(n)$, assign the family 
$((\mathcal{E}^{X}_{j}(\sig))_{j\neq i},\{n\})$
with $\mathcal{E}^{X}_{i}(\sig*(n))=\{n\}$.
For the right upper sequent $\sig*(\ome):\lnot{\rm TI}[<_{{\sf g}_{A}}],\Gam$, 
assign the family $(\mathcal{E}^{X}_{j}(\sig))_{j\neq i}$.
{\small
\[
\infer[(WP)]{(\mathcal{E}^{X}_{j}(\sig))_{j\neq i}\vdash\Gam}
{
\{
((\mathcal{E}^{X}_{j}(\sig))_{j\neq i},\{n\})\vdash\Gam, E_{A}(n)
: n\in\Natural
\}
&
(\mathcal{E}^{X}_{j}(\sig))_{j\neq i}\vdash\lnot{\rm TI}[<_{{\sf g}_{A}}],\Gam
}
\]
}
where
the variable $E_{A}$ with $i=\lceil A\rceil$ does not occur in $\Gam$.

Next suppose that a family 
$((\mathcal{E}^{X}_{j}(\sig))_{j\neq i},\mathcal{E}^{X}_{i}(\sig))\,(i=\lceil A\rceil)$
 is assigned to
the lower sequent $\sig:\Gam,E_{A}(m)$ of the rule $(prg)^{\infty}$.
For each number $n$,
assign the family $((\mathcal{E}^{X}_{j}(\sig))_{j\neq i},\mathcal{E}^{X}_{i}(\sig)\cup\{n\})$
to the $n$-th upper sequent $\sig*(n):\Gam,E_{A}(m),n\not<^{*}_{A}m,E_{A}(n)$ with
$\mathcal{E}^{X}_{i}(\sig*(n))=\mathcal{E}^{X}_{i}(\sig)\cup\{n\}$.
\[
\infer[(prg)^{\infty}]{((\mathcal{E}^{X}_{j}(\sig))_{j\neq i},\mathcal{E}^{X}_{i}(\sig))\vdash\Gam,E_{A}(m)}
{
\{
((\mathcal{E}^{X}_{j}(\sig))_{j\neq i},\mathcal{E}^{X}_{i}\cup\{n\})\vdash
\Gam,E_{A}(m),n\not<^{*}_{A}m,E_{A}(n)
:n\in\Natural
\}
}
\]
where $E_{A}\equiv E_{i}$ with the G\"odel number $i=\lceil A\rceil$ of the formula $A$.

For rules other than $(WP), (prg)^{\infty}$,
the upper sequents receive the same family as the lower sequent receives.
For example
\[
\infer[(\exi^{2}_{1st})]{\mathcal{E}^{X}(\sig)\vdash\exi X F(X),\Gam}
{\mathcal{E}^{X}(\sig)\vdash F(A),\exi X F(X),\Gam}
\]
where $Var(A)\subset Var(\Gam,F)$.

A family $\mathcal{E}^{X}(\sig)$ has been assigned to each node 
$\sig\in T(\pi)$ in the tree of the derivation $\pi$ showing
the fact ${\rm Diag}(\mathcal{P})+(prg)^{\infty}+(WP)\vdash^{\bet}_{\alp}\Gam_{0}$.
}
\edf

Let us define an ordinal function $F(\bet,\alp)$ for giving an upper bound in eliminating
the well-ordering principle.
For normal function ${\sf g}(\alp)$ in Theorems \ref{th:prfthordwop} and \ref{th:derivativemodel}, 
and ordinals $\bet,\alp$,
let us 
define ordinals $F(\bet,\alp)$ recursively on $\alp$ as follows.
$F(\bet,0)=\ome^{1+\bet}$,
\beqn\label{eq:Fsuccessor}
F(\bet,\alp+1)=F\left({\sf g}(\ome^{2(F(\bet,\alp)+\bet)+1})+1+\bet,\alp\right)+{\sf g}(\ome^{2(F(\bet,\alp)+\bet)+1})+1
\eeqn
and $F(\bet,\lam)=\sup\{F(\bet,\alp)+1:\alp<\lam\}$ for limit ordinals $\lam$.

\bprp\label{prp:Fg}
\benu
\item\label{prp:Fg.1}
$\gam<\bet \Rarw F(\gam,\alp)\leq F(\bet,\alp)$, and 
$\gam<\alp\Rarw F(\bet,\gam)<F(\bet,\alp)$.

\item\label{prp:Fg.3}
$F(\bet,\ome(1+\alp))={\sf g}^{\prime}(\alp)$ for $\bet<{\sf g}^{\prime}(\alp)$.

\item\label{prp:Fg.5}
If $\bet<{\sf g}^{\prime}(\alp)$ and $\gam<\ome(1+\alp)$, then $F(\bet,\gam)<{\sf g}^{\prime}(\alp)$.
\eenu
\eprp
\bprf
\ref{prp:Fg}.\ref{prp:Fg.1}.
This follows from the fact that each of functions $\bet\mapsto\alp+\bet$,
$\bet\mapsto\ome^{\bet}$ and $\bet\mapsto{\sf g}(\bet)$
is strictly increasing.
\\
\ref{prp:Fg}.\ref{prp:Fg.5}.
This follows from the fact that ${\sf g}^{\prime}(\alp)$ is closed under $\lam x.\ome^{x}$ and ${\sf g}$.
\eprf
\\

\noindent
The following Elimination theorem \ref{lem:WPbndtransf} of the inference $(WP)$ is a crux for us.

\bth\label{lem:WPbndtransf}{\rm (Elimination of (WP))}\\
Suppose that for a finite set $\Phi$ of $\Sig^{1}_{1}$-formulas
$\vdash^{\bet}_{\alp}\Phi,\Gam$ holds in the calculus
${\rm Diag}(\mathcal{P})+(prg)^{\infty}+(WP)$ 
for $\sig:\Phi,\Gam$ in a witnessing derivation $\pi$
in which the condition (\ref{eq:normal}) is enjoyed for each $(prg)^{\infty}$.
Moreover assume that 
$\Gam\subset\Del(\mathcal{E}^{X}(\sig);\Phi,\Gam)$, and
${\rm Diag}(\mathcal{P},\emptyset)\not\models B$ for any $B\in\Phi$.

Then 
$\mathcal{E}^{X}(\sig)\vdash^{F(\bet,\alp)+\bet}_{0}\Del(\mathcal{E}^{X}(\sig);\Phi,\Gam)$ holds
in the calculus
${\rm Diag}(\mathcal{P},\emptyset)+(prg)^{\emptyset}+(WP)$.
\end{theorem}
In proving Theorem \ref{lem:WPbndtransf},
a key is an extension, Theorem \ref{th:TakeutimultiTI} below,
 of a result due to G. Takeuti\cite{TRemark,PT2}, 
cf.\,Theorem 5 in \cite{Arai17}.

\bth\label{th:TakeutimultiTI}
{\rm The following is provable in ${\rm ACA}_{0}+{\rm WO}(\alp)$:}
\\
Let $n\prec m$ be a binary relation on $\Natural$, 
and $n\prec^{*}m$ the transitive closure of
the relation $n\prec m$.
$(prg)^{D}_{\prec}$ denotes the following inference rule for a predicate $E$.
\[
\infer[(prg)^{D}_{\prec}]{\mathcal{E}^{X}\vdash^{\bet}\Gam,E(m)}
{
\{
\mathcal{E}^{X}\cup\{n\}\vdash^{\bet_{0}}\Gam,E(m),E(n):n\prec^{*}m
\}
}
\]
where the condition (\ref{eq:normal}),
$m\in\mathcal{E}^{X}$, is enjoyed with  a finite set $\mathcal{E}^{X}$.
$\mathcal{E}^{X}\vdash^{\bet}\Gam$ denotes the derivability relation
in a calculus ${\rm Diag}(\emptyset)+(prg)^{D}_{\prec}$,
in which ${\rm Diag}(\emptyset)=\{\bar{E}(n): n\in\Natural\}$.

Assume that there exists an ordinal $\alp$ for which 
$\{n\}\vdash^{\alp}E(n)$ holds
for any natural number $n$.

Then there exist an embedding $f$ such that 
$n\prec m \Rarw f(n)<f(m)$, $f(m)<\ome^{\alp+1}$
for any $n,m\in\Natural$.
\end{theorem}

Proofs of Theorems \ref{lem:WPbndtransf} and \ref{th:TakeutimultiTI}
are postponed in section \ref{sect:CE}.

In what follows we work in ${\rm ACA}_{0}^{+}$.
the set $\{\lceil A\rceil : {\rm Diag}(\mathcal{P},\mathcal{E})\models A\}$ of
the satisfaction relation ${\rm Diag}(\mathcal{P},\mathcal{E})\models A$ for first-order formulas $A$
is then computable from the $\ome$-th jump of the set $\mathcal{P}$.

\subsection{Proof of Theorem \ref{th:prfthordwop}}\label{subsec:Thm3}

First let us prove Theorem \ref{th:prfthordwop}.
In this subsection the predicate $P$ plays no role, and ${\rm Diag}(\mathcal{P})$ is omitted.
Let us introduce a finitary calculus $\bfG_{2}+(prg)+(WPL)$ obtained from
a calculus $\bfG_{2}$ for the predicative second-order logic with inference rules $(\exi^{2}_{1st})$
and $(\fal^{2})$ by adding the following rules $(VJ),(prg),(WPL)$ as follows.
The following inference $(VJ)$ for complete induction schema for first-order formulas $A$
and the successor function $S(x)$ with an eigenvariable $x$.
\[
\infer[(VJ)]{\Gam}
{
\Gam,A(0)
&
\lnot A(x),\Gam,A(S(x))
&
\lnot A(t),\Gam
}
\]
For first-order formulas $A$ and
the eigenvariable $x$:
\[
\infer[(prg)]{\Gam,E_{A}(t)}
{
\Gam,E_{A}(t),x\not<^{*}_{A}t, E_{A}(x)
}
\msfiv
\infer[(WPL)]{\Gam}
{
\Gam, E_{A}(x)
&
\Gam, {\rm LO}(<_{A})
&
\lnot{\rm TI}[<_{{\sf g}_{A}}],\Gam
}
\]
where $E_{A}\equiv E_{i}$ with $i=\lceil A\rceil$, in $(prg)$,
$x$ is the eigenvariable not occurring freely in $\Gam,E_{A}(t)$, and
$Var(A)\subset Var(\Gam)$.
In $(WPL)$, the variable $E_{A}$ does not occur in $\Gam$ (nor in $A$).
The initial sequents are $\Gam,\bar{L},L$ for literals $L$.


We can assume that in a finitary proof,
a variable $E$ occurs in an upper sequent of an inference, but not in the lower sequent
only when the inference is a $(\fal^{2})$, and
the variable $E$ is the eigenvariable of the inference.
Moreover if the end-sequent contains no second-order free variable, then
the variable $E_{A}$ can be assumed not to occur in the lower sequent $\Gam$ of the rule
$(WPL)$ for the relation $<_{A}$.

The axiom of arithmetic comprehension is deduced from the inference rule $(\exi^{2}_{1st})$,
and the axiom ${\rm WOP}({\sf g})$ for the well-ordering principle of ${\sf g}$
is deduced from the inference rules $(prg)$ and  $(WPL)$.

Assume that ${\rm TI}[\prec]$ is provable from ${\rm WOP}({\sf g})$ in ${\rm ACA}_{0}$
for an arithmetical 
relation $\prec$.
Let $\Del_{0}$ denote a set of negations of axioms for first-order arithmetic except complete induction.
By eliminating $(cut)$'s we obtain a proof of $\Del_{0},E_{\prec}(x)$ in 
$\bfG_{2}+(prg)+(WPL)$, where
$E_{\prec}\equiv E_{i}$ with $i=\lceil x_{0}\prec x_{1}\rceil$.


Let us embed the finitary calculus $\bfG_{2}+(prg)+(WPL)$ to 
an intermediate infinitary calculus
$(prg)^{\infty}+(WP)+(cut)_{1st}$,
which is obtained from $(prg)^{\infty}+(WP)$ 
by
adding the cut inference $(cut)_{1st}$ with
a first-order cut formulas $A$:
\[
\infer[(cut)_{1st}]{\Gam,\Del}
{
\Gam,\lnot A^{\infty}
&
A^{\infty},\Del
}
\]

The logical depth $\dg(A)<\ome$ of first-order formulas $A$ is defined recursively by
$\dg(L)=0$ for literals $L$, 
$\dg(A_{0}\lor A_{1})=\dg(A_{0}\land A_{1})=\max\{\dg(A_{0}),\dg(A_{1})\}+1$, and $\dg(\exi x\, A(x)),\dg(\fal x\, A(x))\}=\dg(A(0))+1$.
Then let $\dg(A^{\infty}):=\dg(A)$.
Let $\Gam(x,\ldots)$ be a sequent possibly with free first-order variables $x,\ldots$.
Assuming
$\bfG_{2}+(prg)+(WPL)\vdash\Gam(x,\ldots)$, we see easily
that there exist $d,p,k,m<\ome$ 
such that
$(prg)^{\infty}+(WP)+(cut)_{1st}\vdash^{\ome k+m}_{d,p}\Gam(n,\ldots)^{\infty}$ holds 
for any natural numbers $n,\ldots$,
where the first subscript 
$d$ indicates that the number of nested applications of the rule $(WPL)$
is bounded by $d$, and the second $p$ designates that
any (first-order) cut-formula $A^{\infty}$ occurring in the witnessing derivation has the logical depth
$\dg(A^{\infty})<p$.
 Note that each variable $E$ occurring in the induction formula $A$ of a $(VJ)$ 
 can be assumed to occur
 also in the lower sequent $\Gam$.
We see that there exist $d,p<\ome$ such that
$(prg)^{\infty}+(WP)+(cut)_{1st}\vdash^{\ome^{2}}_{d,p}\Del_{0},E_{\prec}(n)$
holds for any natural number $n$.
Eliminating the false arithmetic $\Del_{0}$,
we obtain
$(prg)^{\infty}+(WP)+(cut)_{1st}\vdash^{\ome^{2}}_{d,p}E_{\prec}(n)$.

Let $2_{0}(\bet)=\bet$ and $2_{p+1}(\bet)=2^{2_{p}(\bet)}$ for $p<\ome$.

\bprp\label{prp:ce1st}
Suppose 
$(prg)^{\infty}+(WP)+(cut)_{1st}\vdash^{\bet}_{d,p}\Gam$.
Then
$(prg)^{\infty}+(WP)\vdash^{2_{p}(\bet)}_{2_{p}(d),0}\Gam$.
\eprp
\bprf
Let $A$ be one of formulas $\exi x\, B$, $B\lor C$, $\bar{E}_{i}(n)$
and arithmetic literals.
We see by induction on $\alp$ that
if $(prg)^{\infty}+(WP)+(cut)_{1st}\vdash^{\bet}_{d,p}\Gam,\lnot A^{\infty}$
and $(prg)^{\infty}+(WP)+(cut)_{1st}\vdash^{\alp}_{e,p}A^{\infty},\Del$
with $\dg(A)\leq p$,
then
$(prg)^{\infty}+(WP)+(cut)_{1st}\vdash^{\bet+\alp}_{d+e,p}\Gam,\Del$.

From the fact we see the proposition by induction on $p<\ome$.
\eprf
\\

\noindent
By Proposition \ref{prp:ce1st} we obtain an ordinal $\bet<\veps_{0}$ and $c<\ome$ for which
$(prg)^{\infty}+(WP)+(cut)_{1st}\vdash^{\bet}_{c,0}E_{\prec}(n)$, i.e.,
$(prg)^{\infty}+(WP)\vdash^{\bet}_{c}E_{\prec}(n)$ holds for any $n$.

In a witnessing derivation $\pi$ of the fact $(prg)^{\infty}+(WP)\vdash^{\bet}_{c}E_{\prec}(n)$,
the condition (\ref{eq:normal}) may be violated.
Let us convert the derivation $\pi$ to a derivation $\pi^{\top}$ as follows.

\bdf
{\rm
For a formula $A$ and a family $\mathcal{E}=(\mathcal{E}_{i})_{i}$ of sets,
$A(\mathcal{E})$ denotes the result of replacing each literal $E_{i}(m)$ by $\top$
when $m\not\in\mathcal{E}_{i}$.
$A(\mathcal{E})$ is defined recursively as follows.
Let $i=\lceil A\rceil$.
\[
\left((E_{A}(m))(\mathcal{E}), (\bar{E}_{A}(m))(\mathcal{E})\right)\equiv\left\{
\begin{array}{ll}
(\top,\bot) & \mbox{if } m\not\in\mathcal{E}_{i}
\\
(E_{A(\mathcal{E})}(m),\bar{E}_{A(\mathcal{E})}(m)) & \mbox{if } m\in\mathcal{E}_{i}
\end{array}
\right.
\]
$L(\mathcal{E})\equiv L$ when $L\in\{\top,\bot,X_{i}(n),\bar{X}_{i}(n):i,n\in\Natural\}$.
$(\bigvee_{n}A_{n})(\mathcal{E}\equiv\bigvee_{n}(A_{n}(\mathcal{E}))$ and
$(\bigwedge_{n}A_{n})(\mathcal{E}\equiv\bigwedge_{n}(A_{n}(\mathcal{E}))$.
$(\exi X\, F(X))(\mathcal{E})\equiv \exi X (F(X)(\mathcal{E}))$ and
$(\fal X\, F(X))(\mathcal{E})\equiv \fal X (F(X)(\mathcal{E}))$.
For a sequent $\Gam$, let $\Gam(\mathcal{E})=\{A(\mathcal{E}):A\in\Gam\}$.

For each node $\sig:\Gam$ in the derivation $\pi$,
let
\[
A^{\sig}:\equiv A(\mathcal{E}^{X}(\sig)), \, \Gam^{\sig}:=\{A^{\sig}:A\in\Gam\}
\]
}
\edf

\bprp\label{prp:top}
Let $\pi$ be a derivation witnessing the fact $\{n\}\vdash^{\bet}_{c}E_{\prec}(n)$ in 
$(prg)^{\infty}+(WP)$.
Then there exists a derivation $\pi^{\top}$ witnessing the same fact in
$(prg)^{\infty}+(WP)$ such that
$\mathcal{E}^{X}_{A^{\sig}}(\sig;\pi^{\top})=\mathcal{E}^{X}_{A}(\sig;\pi)$
for each $\sig\in T(\pi^{\top})\subset T(\pi)$, and
the condition (\ref{eq:normal}) is enjoyed for each rule $(prg)^{\infty}$ occurring in $\pi^{\top}$.
\eprp
\bprf
This is seen by induction on the tree order on the well-founded tree $T(\pi)$.
Each axiom $\sig:\Gam,\bar{L},L$ turns either to $\sig:\Gam^{\sig},\bar{L},L$
or to $\sig:\Gam^{\sig},\bot,\top$.

Consider a rule $(prg)^{\infty}$ in $\pi$.
\[
\infer[(prg)^{\infty}]{((\mathcal{E}^{X}_{j}(\sig;\pi))_{j\neq i},\mathcal{E}^{X}_{i}(\sig;\pi))\vdash\Gam,E_{A}(m)}
{
\{
((\mathcal{E}^{X}_{j}(\sig;\pi))_{j\neq i},\mathcal{E}^{X}_{i}(\sig;\pi)\cup\{n\})\vdash
\Gam,E_{A}(m),n\not<^{*}_{A}m,E_{A}(n)
:n\in\Natural
\}
}
\]
where $E_{A}\equiv E_{i}$ with the G\"odel number $i=\lceil A\rceil$ of the formula $A$.
If $m\in\mathcal{E}^{X}_{A}(\sig;\pi)$, then $m\in\mathcal{E}^{X}_{A^{\sig}}(\sig;\pi^{\top})=\mathcal{E}^{X}_{A}(\sig;\pi)$,
and
\[
\infer[(prg)^{\infty}]{((\mathcal{E}^{X}_{j}(\sig))_{j\neq i},\mathcal{E}^{X}_{i}(\sig))\vdash\Gam^{\sig},E_{A^{\sig}}(m)}
{
\{
((\mathcal{E}^{X}_{j}(\sig))_{j\neq i},\mathcal{E}^{X}_{i}\cup\{n\})\vdash
\Gam^{\sig},E_{A^{\sig}}(m),n\not<^{*}_{A^{\sig}}m,E_{A^{\sig}}(n)
:n\in\Natural
\}
}
\]
Otherwise $(E_{A}(m))^{\sig}\equiv\top$.
$\Gam^{\sig},\top$ is an axiom. Discard the upper part.

From the construction of $\pi^{\top}$ we see easily that
$\mathcal{E}^{X}_{A^{\sig}}(\sig;\pi^{\top})=\mathcal{E}^{X}_{A}(\sig;\pi)$ for each node 
$\sig\in T(\pi^{\top})\subset T(\pi)$, and
the condition (\ref{eq:normal}) is enjoyed for each rule $(prg)^{\infty}$ occurring in $\pi^{\top}$.
\eprf
\\

By Proposition \ref{prp:top}
we obtain a derivation $\pi^{\top}$ witnessing the fact $\{n\}\vdash^{\bet}_{c}E_{\prec}(n)$ in
$(prg)^{\infty}+(WP)$ such that
the condition (\ref{eq:normal}) is enjoyed for each rule $(prg)^{\infty}$ occurring in $\pi^{\top}$.

We see from Theorem \ref{lem:WPbndtransf}
 that in the calculus ${\rm Diag}(\emptyset)+(prg)^{\emptyset}$,
$\{n\}\vdash^{\alp}_{0}E_{\prec}(n)$ holds for any $n$, and the ordinal
$\alp=F(\bet,c)+\bet$, where 
$\{n\}=\mathcal{E}_{i}^{X}(\pi)$ with $i=\lceil x_{0}\prec x_{1}\rceil$ and
$\Del(\{n\},\emptyset;E_{\prec}(n))=\{E_{\prec}(n)\}$.

Theorem \ref{th:TakeutimultiTI} yields an embedding
$f$ such that $n\prec m \Rarw f(n)<f(m)<\ome^{\alp+1}$.

On the other hand we have $\ome^{\alp+1}=\ome^{F(\bet,c)+\bet+1}<{\sf g}^{\prime}(0)$
by Proposition \ref{prp:Fg}.\ref{prp:Fg.5} and $\bet<\veps_{0}\leq{\sf g}^{\prime}(0)$.
Thus Theorem \ref{th:prfthordwop} is proved.

\subsection{Corrections to \cite{Arai17}}\label{sect:correction}
The proof of the harder direction of Theorem 4 in \cite{Arai17} 
should be corrected as pointed out by A. Freund.
In this subsection the predicate $P$ will denote a given set of natural numbers.
Let us augment another countable list $Y_{i},\bar{Y}_{i}\,(i\in\Natural)$ of
second-order free variables.
First-order formulas may contain these variables $Y$.

Assuming ${\rm WOP}({\sf g}^{\prime})$, we need to show the existence of a countable coded $\ome$-model 
$(\mathcal{P},(\mathcal{Q})_{i})_{i<\ome}$ of ${\rm ACA}_{0}+{\rm WOP}({\sf g})$ for a given set $\mathcal{P}\subset\Natural$.
In what follows argue in ${\rm ACA}_{0}+{\rm WOP}({\sf g}^{\prime})$.
Since ${\rm WOP}({\sf g}^{\prime})$ implies ${\rm WOP}(\lam X.\veps_{X})$, which in turn
yields ${\rm ACA}_{0}^{+}$ by Theorem \ref{th:MontalbanMarcone},
we are working in ${\rm ACA}_{0}^{+}+{\rm WOP}({\sf g}^{\prime})$, and we can assume the
existence of the $\ome$-th jump of any sets.

Let us search a derivation of the contradiction $\emptyset$ in the following infinitary calculus
${\rm Diag}(\mathcal{P})+(prg)^{\infty}+(WP)+(ACA)$, in which the variables $E_{i}$ are not interpreted.
The calculus is obtained from the infinitary 
calculus $(prg)^{\infty}+(WP)$
 by
 adding the following rule $(ACA)$ for arithmetic comprehension axiom:
\[
\infer[(ACA)]{\Gam}
{
Y_{j}\neq A,\Gam
}
\]
where $A$ is a first-order formula, 
$Y_{j}$ is the eigenvariable not occurring in $\Gam\cup\{A\}$, and
$Y_{j}\neq A:\Lrarw (\lnot\fal x[Y_{j}(x)\lrarw A(x)])^{\infty}$.
Note that variables $E_{i},Y_{i}$ are uninterpreted in the calculus.

A tree $\calt\subset{}^{<\ome}\Natural$ is constructed recursively as follows.
At each node $\sig$, a sequent and a family $\mathcal{E}^{X}(\sig)$ of finite sets are assigned.
At the bottom $\emptyset$, we put the empty sequent, and $\mathcal{E}^{X}(\sig)=\emptyset$.
The assignment $\mathcal{E}^{X}(\sig)$ is done similarly as in Definition \ref{df:attach}.

Suppose that the tree $\calt$ has been constructed up to a node $\sig\in{}^{<\ome}\Natural$. 
Let $\{A_{i}\}_{i}$ be an enumeration of all first-order formulas (abstracts).
\\
{\bf Case 0}. The length $lh(\sig)=3i$:
Apply one of inferences $(\bigvee),(\bigwedge),(\exi^{2}_{1st})$,
and $(prg)^{\infty}$ if it is possible.
Otherwise repeat, i.e., apply an inference $(Rep)$.

When $(\exi^{2}_{1st})$ is applied backwards, a first-order $A_{j}$ is chosen so that $j$ is the least 
such that
$A_{j}$ has not yet been tested for the major formula $\exi X\, F(X)$ of the $(\exi^{2}_{1st})$, and
$Var(A_{j})\subset Var(\Gam\cup\{F\})\cup\{P\}$.
\[
\infer[(\exi^{2}_{1st})]{\Gam,\exi X\, F(X)}
{
\Gam,\exi X\,F(X),F(A_{j})
}
\]
When $(prg)^{\infty}$ is applied backwards to a formula $E_{A}(m)$ with $i=\lceil A\rceil$,
the condition (\ref{eq:normal}), $m\in\mathcal{E}^{X}_{i}(\sig)$,
and $Var(A)\subset Var(\Gam)$
 have to be met.
Otherwise repeat.
\[
\infer[(prg)^{\infty}]{((\mathcal{E}^{X}_{j}(\sig))_{j\neq i},\mathcal{E}^{X}_{i}(\sig))\vdash\Gam,E_{A}(m)}
{
\{
((\mathcal{E}^{X}_{j}(\sig))_{j\neq i},\mathcal{E}^{X}_{i}(\sig)\cup\{n\})\vdash\Gam,n\not<^{*}_{A}m, E_{A}(n)
:n\in\Natural
\}
}
\]
{\bf Case 1}. $lh(\sig)=3\la i,n\ra+1$: 
Apply the inference $(ACA)$ backwards with the first-order $A\equiv A_{i}$
and an eigenvariable $Y_{j}$
 if $Var(A)\subset Var(\Gam)\cup\{P\}$.
 \[
\infer[(ACA)]{\Gam}
{
Y_{j}\neq A,\Gam
}
\]
Otherwise repeat
\\
{\bf Case 2}. $lh(\sig)=3i+2$:
Apply the inference $(WP)$ backwards with the relation $<_{A_{i}}$.
\\

If the tree $\calt$ is not well-founded, then let $\mathcal{R}$ be an infinite path through $\calt$.
We see for any $i,n$ that
at most one of 
 $Q(n)$ or $\bar{Q}(n)$ is on $\mathcal{R}$
 for $Q\in\{E_{i},Y_{i},P:i\in\Natural\}$,
and 
$[(P(n))\in\mathcal{R}\Rarw n\not\in\mathcal{P}]\spand[(\bar{P}(n))\in\mathcal{R}\Rarw n\in\mathcal{P}]$
 due to the axioms
$\Gam,L$ with $L\in{\rm Diag}(\mathcal{P})$.
Let $(\mathcal{Q})_{i}$ be the set defined by
$n\in(\mathcal{Q})_{2i} \Lrarw (E_{i}(n))\not\in\mathcal{R}$
and
$n\in(\mathcal{Q})_{2i+1} \Lrarw (Y_{i}(n))\not\in\mathcal{R}$.

$(\mathcal{P},(\mathcal{Q})_{i})_{i\in\Natural}$ is shown to be a countable coded $\ome$-model of 
${\rm ACA}_{0}+{\rm WOP}({\sf g})$ as follows.
The search procedure is fair, i.e., each formula is eventually analyzed on every path.
To ensure fairness, formulas in sequents $\Gam$ are assumed to stand in a queue.
The head of the queue is analyzed in {\bf Case 0}, and the analyzed formula moves to the end of the queue
in the next stage.
We see from the fairness that
${\rm Diag}(\mathcal{P},(\mathcal{Q})_{i})_{i\in\Natural}\not\models A$ 
first by induction on the number of occurrences of 
logical connectives in first-order
formulas $A$ on the path $\mathcal{R}$,
and then for $\Pi^{1}_{1}$-formulas ${\rm TI}[<_{A}]$ and $\Sig^{1}_{1}$-formulas
$\lnot{\rm TI}[<_{{\sf g}_{A}}]$.
Moreover
${\rm Diag}(\mathcal{P},(\mathcal{Q})_{i})_{i\in\Natural}\models{\rm ACA}_{0}$ since the inference rules $(ACA)$
are analyzed for every $A_{i}$.
Finally we show
${\rm Diag}(\mathcal{P},(\mathcal{Q})_{i})_{i\in\Natural}\models{\rm WOP}({\sf g})$.
Assume that ${\rm Diag}(\mathcal{P},(\mathcal{Q})_{i})_{i\in\Natural}\models{\rm WO}[<_{A}]$
for a first-order $A$.
The path $\mathcal{R}$ passes through an inference $(WP)$ for the relation $<_{A}$.
If $\mathcal{R}$ passes through the rightmost upper sequent
$\lnot{\rm TI}[<_{{\sf g}_{A}}]$, then
 ${\rm Diag}(\mathcal{P},(\mathcal{Q})_{i})_{i\in\Natural}\not\models\lnot{\rm TI}[<_{{\sf g}_{A}}]$,
 i,e.,  ${\rm Diag}(\mathcal{P},(\mathcal{Q})_{i})_{i\in\Natural}\models{\rm TI}[<_{{\sf g}_{A}}]$,
 and we are done.
Suppose that $\mathcal{R}$ passes through an $n_{0}$-th upper sequent
$((\mathcal{E}_{j}^{X}(\sig_{0}))_{j\neq i},\mathcal{E}^{X}_{i}(\sig_{0}))\vdash\Gam_{0},E_{A}(n_{0})$ 
and $E_{A}=E_{i}$ with $i=\lceil A\rceil$.
Since the condition (\ref{eq:normal}), $n_{0}\in\mathcal{E}^{X}_{i}(\sig_{0})=\{n_{0}\}$ is met,
the formula $E_{A}(n_{0})$ is analyzed after a number of steps at a $(prg)^{\infty}$,
and $\mathcal{R}$ passes through an $n_{1}$-th branch
$((\mathcal{E}_{j}^{X}(\sig_{1}))_{j\neq i},\mathcal{E}^{X}_{i}(\sig_{1}))\vdash\Gam_{1},
n_{1}\not<_{A}^{*}n_{0},E_{A}(n_{1}),E_{A}(n_{0})$.
We obtain ${\rm Diag}(\mathcal{P},(\mathcal{Q})_{i})_{i\in\Natural}\not\models n_{1}\not<_{A}^{*}n_{0}$,
 i.e., $n_{1}<_{A}^{*,\mathcal{P},\mathcal{Q}}n_{0}$.
 Also $\{n_{0},n_{1}\}\subset\mathcal{E}_{i}^{X}(\sig_{1})$.
In this way we obtain an infinite descending chain 
$\cdots<_{A}^{*,\mathcal{P},\mathcal{Q}}n_{2}<_{A}^{*,\mathcal{P},\mathcal{Q}}n_{1}<_{A}^{*,\mathcal{P},\mathcal{Q}}n_{0}$ from $\mathcal{R}$, contradicting the assumption 
${\rm WO}[<^{\mathcal{P},\mathcal{Q}}]$.
\\

In what follows assume that the tree $\calt$ is well-founded.
Let $\Lam$ denote the least epsilon number
larger than the order type of the Kleene-Brouwer ordering $<_{KB}$
on the well-founded tree $\calt$.
We have ${\rm WO}({\sf g}^{\prime}(\Lam))$
by ${\rm WOP}({\sf g}^{\prime})$ and ${\rm WO}(\Lam)$.

For $b<\Lam$  let us write $S+(ACA)\vdash^{b}_{c}\Gam$ 
when there exists a derivation of $\Gam$ in 
${\rm Diag}(\mathcal{P})+(prg)^{\infty}+(WP)+(ACA)$ such that
its depth is bounded by $b$, the depth of nested applications of the rules $(WP)$
 is bounded by $c$, and the condition (\ref{eq:normal}) is enjoyed for each 
 inference $(prg)^{\infty}$ in the derivation, where a family $\mathcal{E}^{X}(\sig)$ of finite sets
 is assigned to each node $\sig$ in the derivation tree as in Definition \ref{df:attach}.
 
For the inference
\[
\infer[(ACA)]{\Gam}
{
Y_{j}\neq A,\Gam
}
\]
substitute $A$ for the eigenvaraible $Y_{j}$, 
and deduce the valid formula $A= A$ logically in a finite number of steps, and
then a $(cut)_{1st}$ yields the lower sequent $\Gam$.
Axioms $\Gam,\bar{Y}_{j}(n),Y_{j}(n)$ turns to another valid sequent
$\Gam^{\prime}, \lnot A(n),A(n)$.
In $(WP)$, if $Y_{j}$ occurs in $B(Y_{j})$, then the variable $E_{B(Y_{j})}\equiv E_{i}$ with 
$i=\lceil B(Y_{j})\rceil$
should be renamed to $E_{B(A)}\equiv E_{k}$ with $k=\lceil B(A)\rceil$.
{\small
\[
\infer{\Gam(Y_{j})}
{
\Gam(Y_{j}), E_{B(Y_{j})}(n)
&
\lnot{\rm TI}[<_{{\sf g}_{B(Y_{j})}}],\Gam(Y_{j})
}
\leadsto
\infer{\Gam(A)}
{
\Gam(A), E_{B(A)}(n)
&
\lnot{\rm TI}[<_{{\sf g}_{B(A)}}],\Gam(A)
}
\]
}
Thus we obtain 
${\rm Diag}(\mathcal{P})+(prg)^{\infty}+(WP)+(cut)_{1st}\vdash^{\ome+b}_{b}\emptyset$
from ${\rm Diag}(\mathcal{P})+(prg)^{\infty}+(WP)+(ACA)\vdash^{b}_{b}\emptyset$,
and
${\rm Diag}(\mathcal{P})+(prg)^{\infty}+(WP)\vdash^{2_{p}(\ome+b)}_{2_{p}(b)}\emptyset$
for a $p<\ome$ as in Proposition \ref{prp:ce1st}.

Here the condition (\ref{eq:normal}) is forced in the search.

Theorem \ref{lem:WPbndtransf} yields 
$\emptyset\vdash^{\del}_{0}\emptyset$
with $\del=F(2_{p}(\ome+b),2_{p}(b))+2_{p}(\ome+b)$
 and $\mathcal{E}^{X}(\emptyset)=\Del(\emptyset;\emptyset)=\emptyset$.
This means that in the $\ome$-logic, there exists a cut-free derivation of $\emptyset$ in depth
$\del<{\sf g}^{\prime}(\Lam)$, which is seen 
 from Proposition \ref{prp:Fg}.\ref{prp:Fg.5} and $b<\Lam$.
We see by induction up to the ordinal ${\sf g}^{\prime}(\Lam)$ that
this is not the case.
Therefore the tree $\calt$ must not be well-founded.
Thus our proof of Theorem \ref{th:derivativemodel} is completed.

\section{Elimination of the inference for well-ordering principle}\label{sect:CE}
It remains to show Theorems \ref{th:TakeutimultiTI} and \ref{lem:WPbndtransf}.
\\

\noindent
({\bf Proof} of Theorem \ref{th:TakeutimultiTI}).
We can assume that the transitive closure $\prec^{*}$ of the relation $n\prec m$
is irreflexive.
Namely there is no sequence $(n_{0},\ldots,n_{k})\,(k\geq 1)$
such that $n=n_{0}=n_{k}$ and $\fal j<k(n_{j+1}\prec n_{j})$.
Suppose that there exists such a sequence.
By the assumption we obtain
$\{n_{0}\}\vdash^{\alp_{0}}E(n_{0})$ for $\alp_{0}=\alp$.
Any positive literal $E(n_{0})$ is not an axiom in 
${\rm Diag}(\emptyset)$.
We see from $n_{0}\in\{n_{0}\}$ by induction on ordinals $\alp_{0}$ that
there must be an inference $(prg)_{\prec}^{D}$ in the witnessed derivation,
and we obtain $\{n_{0},n_{1}\}\vdash^{\alp_{1}}E(n_{0}),E(n_{1})$ for 
an $\alp_{1}<\alp_{0}$.
Again $P(n_{0}),P(n_{1})$ is not an axiom in
${\rm Diag}(\emptyset)$.
In this way we would obtain an infinite descending chain $\{\alp_{m}\}_{m<\ome}$
of ordinals such that  $\{n_{j}\}_{j<k}\vdash^{\alp_{m}}_{0}\{P(n_{j}):j<k\}$.



By recursion on $m$, we define 
a non-empty finite set $\mathcal{E}(m)$, and
an ordinal $\bet(m)\leq\alp$ for which the followings hold for
$\Del(\mathcal{E}(m)):=\{E(n): n\in\mathcal{E}(m)\}$.
\beqnarr
&&
\mathcal{E}(m)\subset\{n: m\preceq^{*}n\leq m\}
\spand
\mathcal{E}(m)\vdash^{\bet(m)}\Del(\mathcal{E}(m))
\nonumber
\\
&&
\fal n<m(m\prec^{*}n \to \bet(m)<\bet(n))
\label{eq:Takeutimulti}
\eeqnarr
{\bf Case 1}. $\lnot\exi n<m(m\prec^{*}n)$:
Let $\mathcal{E}(m)=\{m\}$
and $\bet(m)=\alp$.
Then the conditions in (\ref{eq:Takeutimulti}) are fulfilled with 
$\Del(\mathcal{E}(m))=\{E(m)\}$.
\\
{\bf Case 2}. Otherwise:
Pick a $k<m$ such that $m\prec^{*}k$ and
$\bet(k)=\min\{\bet(n): n<m, m\prec^{*}n\}$.

Then let $\mathcal{E}(m)=\mathcal{E}(k)\cup\{m\}$.
On the other hand we have $\mathcal{E}(k)\vdash^{\bet(k)}\Del(\mathcal{E}(k))$.
The sequent $\Del(\mathcal{E}(k))$ is not an axiom in ${\rm Diag}(\emptyset)$.
Search the lowest inference $(prg)_{\prec}^{D}$ 
in the derivation showing the fact $\mathcal{E}(k)\vdash^{\bet(k)}\Del(\mathcal{E}(k))$:
\[
 \infer[(prg)_{\prec}^{D}]{\mathcal{E}(k)\vdash^{\bet^{\prime}}\Del(\mathcal{E}(k))}
 {
  \{
  \mathcal{E}(k)\cup\{n\}\vdash^{\bet_{0}}\Del(\mathcal{E}(k)), E(n): n\prec^{*} k^{\prime}
   \}
}
\]
where 
$\bet_{0}<\bet^{\prime}\leq\bet(k)$,
there may be some $(Rep)$'s below the inference $(prg)_{\prec}^{D}$,
and
$E(k^{\prime})\in\Del(\mathcal{E}(k))$ is the main formula of the inference 
$(prg)_{\prec}^{D}$.
We have $m\prec^{*}k\preceq^{*}k^{\prime}$, and $m\prec^{*}k^{\prime}$.
Pick the $m$-th branch in the upper sequents.
We obtain $\mathcal{E}(m)\vdash^{\bet(m)}\Del(\mathcal{E}(m))$
for $\bet(m):=\bet_{0}<\bet(k)$.
The conditions in (\ref{eq:Takeutimulti}) are fulfilled.

Now define a function $f(m)$ as follows.
\[
f(m)=\max\{\ome^{\bet(m_{0})}\#\cdots\#\ome^{\bet(m_{k})}: m_{0}\prec^{*}\cdots\prec^{*}m_{k}=m, 
m_{0},\ldots,m_{k-1}<m\}
\]
where $\#$ denotes the natural sum.
Note that the set 
$\{(m_{0},\ldots,m_{k}): m_{0}\prec^{*}\cdots\prec^{*}m_{k}=m, 
m_{0},\ldots,m_{k-1}<m\}$ is finite since $\prec^{*}$ is irreflexive.

We show the function $f$ is a desired embedding between $\prec^{*}$ and $<$.
Assume $m\prec^{*}n$, and let $m_{0},\ldots,m_{k}$ be a sequence such that
$f(m)=\ome^{\bet(m_{0})}\#\cdots\#\ome^{\bet(m_{k})}$,
with $m_{0}\prec^{*}\cdots\prec^{*}m_{k}=m$ and 
$m_{0},\ldots,m_{k-1}<m$.
We obtain $m_{i}\prec^{*}n$ for any $i\leq k$.
Let us partition the set $\{0,\ldots,k\}$ into two sets
$A=\{i\leq k: n<m_{i}\}$ and $B=\{i\leq k: m_{i}<n\}$.
Note that $m_{i}\neq n$ since $\prec^{*}$ is irreflexive.

By (\ref{eq:Takeutimulti}) we obtain $\bet(m_{i})<\bet(n)$ for each $i\in A$, and hence
$\#\{\ome^{\bet(m_{i})}:i\in A\}<\ome^{\bet(n)}$,
where
$\#\{\alp_{1},\ldots,\alp_{n}\}=\alp_{1}\#\cdots\#\alp_{n}$.

On the other hand we have 
$\#\{\ome^{\bet(m_{i})}:i\in B\}\leq
\max\{\ome^{\bet(n_{0})}\#\cdots\#\ome^{\bet(n_{\ell})}: n_{0}\prec^{*}\cdots\prec^{*}n_{\ell-1}, 
n_{0},\ldots,n_{\ell-1}<n\}$.
Therefore we conclude
$f(m)<f(n)$.
\eprf

\blem\label{lem:bndwset}
For each $j\leq\ell$, let $<_{j}$ be a first-order formula with $j=\lceil <_{j}\rceil$.
Let $\mathcal{E}^{X}=(\mathcal{E}^{X}_{j})_{j<\ell}$ be finite
sets, and $\mathcal{E}^{X}_{\ell}$ a finite set.
Let $\Gam\subset\bigcup_{j<\ell}\Del(\mathcal{E}^{X}_{j})$ be a sequent
and
$\Gam_{\ell}\subset\Del(\mathcal{E}^{X}_{\ell})=\{E_{\ell}(n): n\in\mathcal{E}^{X}_{\ell}\}$
a sequent.
In the calculus ${\rm Diag}(\mathcal{P},\emptyset)+(prg)^{\mathcal{E}}$,
assume that 
$(\mathcal{E}^{X},\mathcal{E}^{X}_{\ell})\vdash^{\alp}\Gam,\Gam_{\ell}$
for $\mathcal{E}^{X}=(\mathcal{E}^{X}_{j})_{j<\ell}$.
Then either
$\mathcal{E}^{X}\vdash^{\alp}\Gam$ holds, or
$\mathcal{E}^{X}_{\ell}\vdash^{\alp}\Gam_{\ell}$ holds.
\elem
\bprf
We show the lemma by induction on $\alp$.
Assume that
$(\mathcal{E}^{X},\mathcal{E}^{X}_{\ell})\vdash^{\alp}\Gam$ does not hold.
The set $\Gam\cup\Gam_{\ell}$ consisting of positive literals $E_{i}(n)$, is not an axiom
in ${\rm Diag}(\mathcal{P},\emptyset)+(prg)^{\emptyset}$.

Consider the case when the last inference is a $(prg)^{\emptyset}$
for a $<_{j}$:
\[
\infer[(prg)^{\emptyset}]{(\mathcal{E}^{X},\mathcal{E}^{X}_{\ell})\vdash^{\alp}\Gam,\Gam_{\ell}}
{
\{
(\mathcal{E}^{X},\mathcal{E}^{X}_{\ell})_{j,m}\vdash^{\bet}\Gam,\Gam_{\ell},E_{j}(\bar{m})
: m<^{*,\emptyset}_{j} n
\}
}
\]
where $E_{j}(\bar{n})$ is in $\Gam\cup\Gam_{\ell}$,
and $(\mathcal{E}^{X},\mathcal{E}^{X}_{\ell})_{j,m}$ denotes
the sequence $(\mathcal{E}^{X},\mathcal{E}^{X}_{\ell})$ except
$\mathcal{E}^{X}_{j}$ is replaced by $\mathcal{E}^{X}_{j}\cup\{m\}$.
We have (\ref{eq:normal}), $n\in\mathcal{E}^{X}_{j}$.

First consider the case $j\neq\ell$.
By the assumption we see that there exists an 
$m<^{*,\emptyset}_{j} n$
such that $(\mathcal{E}^{X},\mathcal{E}^{X}_{\ell})_{j,m}\vdash^{\bet}\Gam,E_{j}(\bar{m}),E_{j}(\bar{n})$
does not hold. IH yields 
$(\mathcal{E}^{X},\mathcal{E}^{X}_{\ell})\vdash^{\bet}\Gam_{\ell}$.

Second consider the case $j=\ell$.
We see from the assumption that
for each $m<^{*,\emptyset}_{\ell} n$,
$(\mathcal{E}^{X},\mathcal{E}^{X}_{\ell})_{\ell,m}\vdash^{\bet}\Gam_{\ell},E_{\ell}(\bar{m})$ holds.
Then an inference $(prg)^{\emptyset}$ yields
$(\mathcal{E}^{X},\mathcal{E}^{X}_{\ell})\vdash^{\alp}\Gam_{\ell}$.
\eprf
\\

\noindent
({\bf Proof} of Theorem \ref{lem:WPbndtransf}).
Let us prove Theorem \ref{lem:WPbndtransf}
by induction on $\bet$.
Suppose that in the calculus ${\rm Diag}(\mathcal{P})+(prg)^{\infty}+(WP)$,
$\vdash^{\bet}_{\alp}\Phi,\Gam$ holds in a derivation $\pi$
for a sequent $\sig:\Phi,\Gam$ such that
$\Phi$ is a $\Sig^{1}_{1}$-formulas with 
${\rm Diag}(\mathcal{P},\emptyset)\not\models \bigvee\Phi$ and
$\Gam\subset\Del(\mathcal{E}^{X}(\sig);\Phi,\Gam)$.
Moreover the condition (\ref{eq:normal}) is assumed to be
enjoyed for each $(prg)^{\infty}$
occurring in $\pi$.
Let $\mathcal{E}^{X}=\mathcal{E}^{X}(\sig)$.

If $\Phi,\Gam$ is an axiom in ${\rm Diag}(\mathcal{P})+(prg)^{\infty}+(WP)$,
then $\Del(\mathcal{E}^{X};\Gam)$ is an axiom in ${\rm Diag}(\mathcal{P},\emptyset)+(prg)^{\mathcal{E}}$
since ${\rm Diag}(\mathcal{P},\emptyset)\not\models\bigvee\Phi$.
For example, consider the case when $\{\bar{E}_{i}(m),E_{i}(m)\}\subset\Phi\cup\Gam$.
Since $\Gam$ contains positive literals only,
we may assume $\bar{E}_{i}(m)\in\Phi$ and $E_{i}(m)\in\Gam$.
From ${\rm Diag}(\mathcal{P},\emptyset)\models\bar{E}_{i}(m)$ 
we see that this is not the case.

Consider the last inference in the derivation showing 
$\vdash^{\bet}_{\alp}\Phi,\Gam$.
\\
{\bf Case 1}. The last inference is a $(prg)^{\infty}$. For $\gam<\bet$
and $i=\lceil A\rceil$, we have for an $E_{i}(m)\in\Phi\cup\Gam$
\[
\infer[(prg)^{\infty}]{((\mathcal{E}^{X}_{j})_{j\neq i},\mathcal{E}^{X}_{i})\vdash^{\bet}_{\alp}\Phi,\Gam}
{
\{
((\mathcal{E}^{X}_{j})_{j\neq i},\mathcal{E}^{X}_{i}\cup\{n\})\vdash^{\gam}_{\alp}
\Phi,n\not<^{*}_{A}m,\Gam,E_{i}(n):
n\in\Natural
\}
}
\]
IH yields 
$((\mathcal{E}^{X}_{j})_{j\neq i},\mathcal{E}^{X}_{i}\cup\{n\})\vdash^{F(\gam,\alp)+\gam}
\Del_{n}$
for each $n<^{*,\emptyset}_{A}m$,
where $E_{i}(n)\in \Del_{n}=\Del((\mathcal{E}^{X}_{j})_{j\neq i},\mathcal{E}^{X}_{i}\cup\{n\};\Phi,n\not<^{*}_{A}m,\Gam,E_{i}(m),E_{i}(n))$,
and 
$\Del((\mathcal{E}^{X}_{j})_{j\neq i},\mathcal{E}^{X}_{i}\cup\{n\};\Phi,\Gam,E_{i}(m))\cup\{E_{i}(n)\}=\Del_{n}
$.
By the assumption we obtain (\ref{eq:normal}), $m\in\mathcal{E}^{X}_{i}$.
An inference $(prg)^{\emptyset}$ with $F(\gam,\alp)+\gam<F(\bet,\alp)+\bet$ 
yields
$((\mathcal{E}^{X}_{j})_{j\neq i},\mathcal{E}^{X}_{i})\vdash^{F(\bet,\alp)+\bet}
\Del((\mathcal{E}^{X}_{j})_{j\neq i},\mathcal{E}^{X}_{i};\Phi,\Gam,E_{i}(m))$.
\\

\noindent
{\bf Case 2}. The last inference is a $(WP)$. 
For $\gam<\bet$, $\alp_{0}<\alp$, $Var(A)\subset Var(\Gam)\cup\{P\}$, $i=\lceil A\rceil$, and
the variable $E_{i}$ not occurring in $\Gam$, we have
\[
\infer[(WP)]{(\mathcal{E}^{X}_{j})_{j\neq i}\vdash^{\bet}_{\alp}\Phi,\Gam}
{
\{
((\mathcal{E}^{X}_{j})_{j\neq i},\{n\})\vdash^{\gam}_{\alp_{0}}\Phi,\Gam, E_{i}(n)
: n\in\Natural
\}
&
\mathcal{E}^{X}\vdash^{\gam}_{\alp_{0}}\lnot{\rm TI}[<_{{\sf g}_{A}}],\Phi,\Gam
}
\]

For the left upper sequent we have for each $n\in\Natural$
\[
((\mathcal{E}^{X}_{j})_{j\neq i},\{n\})\vdash^{\gam}_{\alp_{0}}\Phi,\Gam, E_{i}(n)
\]
By IH we obtain 
$((\mathcal{E}^{X}_{j})_{j\neq i},\{n\})\vdash^{F(\gam,\alp_{0})+\gam}\Del((\mathcal{E}^{X}_{j})_{j\neq i};\Phi,\Gam),E_{i}(n)$.
Here note that
 $\Del((\mathcal{E}^{X}_{j})_{j\neq i};\Phi,\Gam)\cup\{E_{i}(n)\}=\Del((\mathcal{E}^{X}_{j})_{j\neq i},\{n\};\Phi,\Gam,E_{i}(n))$.
From Lemma \ref{lem:bndwset} we see that
either
$(\mathcal{E}^{X}_{j})_{j\neq i}\vdash^{F(\gam,\alp_{0})+\gam}\Del(\mathcal{E}^{X};\Phi,\Gam)$
or
$((\mathcal{E}^{X}_{j})_{j\neq i},\{n\})\vdash^{F(\gam,\alp_{0})+\gam}E_{i}(n)$ holds.
If $(\mathcal{E}^{X}_{j})_{j\neq i}\vdash^{F(\gam,\alp_{0})+\gam}\Del(\mathcal{E}^{X};\Phi,\Gam)$, then we are done.

Assume $((\mathcal{E}^{X}_{j})_{j\neq i},\{n\})\vdash^{F(\gam,\alp_{0})+\gam} E_{i}(n)$,
i.e., $\{n\}\vdash^{F(\gam,\alp_{0})+\gam} E_{i}(n)$
 for every $n$.
By Theorem \ref{th:TakeutimultiTI} we obtain
an embedding $f$
from $<_{A}^{\emptyset}$ to $\ome^{F(\gam,\alp_{0})+\gam+1}$, which yields
an embedding from $<_{{\sf g}(A)}^{\emptyset}$
to $\del:={\sf g}(\ome^{F(\gam,\alp_{0})+\gam+1})$ by Proposition \ref{prp:dilate}.
Hence we see from ${\rm WO}(\del)$ that
${\rm Diag}(\mathcal{P},\emptyset)\models {\rm TI}[<_{{\sf g}_{A}},C]$
for any first-order formulas $C$.
Therefore ${\rm Diag}(\mathcal{P},\emptyset)\models \lnot{\rm TI}[<_{{\sf g}_{A}}]$.

Second consider the right upper sequent. 
IH yields 
\[
(\mathcal{E}^{X}_{j})_{j\neq i}\vdash^{F(\gam,\alp_{0})+\gam}\Del(\mathcal{E}^{X};\Phi,\Gam)
\]
On the other side Proposition \ref{prp:Fg}.\ref{prp:Fg.1} with (\ref{eq:Fsuccessor}) yields 
$F(\gam,\alp_{0})+\gam\leq F(\gam,\alp)+\gam\leq F(\bet,\alp)+\bet$.
Hence the assertion $\mathcal{E}^{X}\vdash^{F(\bet,\alp)+\bet}\Del(\mathcal{E}^{X};\Phi,\Gam)$ follows.
\\
{\bf Case 3}. The last inference is other than $(prg)^{\infty}$ and $(WP)$.

Consider first the case when the last inference is a rule for existential second-order quantifier.
\[
\infer[(\exi_{1st}^{2})]{(\mathcal{E}^{X}_{j})_{j\neq i}\vdash^{\bet}_{\alp}\lnot{\rm TI}[<_{{\sf g}_{A}}],\Phi,\Gam}
{
(\mathcal{E}^{X}_{j})_{j\neq i}\vdash^{\gam}_{\alp} \lnot{\rm TI}[<_{{\sf g}_{A}}],\lnot {\rm TI}[<_{{\sf g}_{A}},C],\Phi,\Gam
}
\]
where $\gam<\bet$ and $C$ is a first-order formula such that
$Var(C)\subset Var(A,\Gam)$.
IH yields
$(\mathcal{E}^{X}_{j})_{j\neq i}\vdash^{F(\gam,\alp)+\gam} \Del(\mathcal{E}^{X};\Phi,\Gam)$ for
$\Del(\mathcal{E}^{X};\Phi,\Gam)=\Del(\mathcal{E}^{X};\lnot {\rm TI}[<_{{\sf g}_{A}},C],\Phi,\Gam)$
by $Var(A,C)\subset Var(\Phi,\Gam)$.

Second consider the case when the last inference is a rule $(\bigwedge)$
For $\gam<\bet$ we have
\[
\infer[(\bigwedge)]{\mathcal{E}^{X}\vdash^{\bet}_{\alp}\Phi,\Gam}
{
\{
\mathcal{E}^{X}\vdash^{\gam}_{\alp}A_{n},\Phi,\Gam
: n\in\Natural
\}
}
\]
where $\bigwedge_{n}A_{n}$ is in the set $\Phi$.
Let $n$ be the least number such that ${\rm Diag}(\mathcal{P},\emptyset)\not\models A_{n}$.
IH yields 
$\mathcal{E}^{X}\vdash^{F(\gam,\alp)+\gam}\Del(\mathcal{E}^{X};\Phi,\Gam)$.

Other cases are similarly seen.
\eprf

\end{document}